\documentclass{llncs}
\usepackage[centertags]{amsmath}
\usepackage{amssymb}
\usepackage[dvipdfm,bookmarks=true,colorlinks=true,linkcolor=blue,urlcolor=green,citecolor=red,pdfstartview=FitH]{hyperref} %,pagebackref

\begin{document}

\title{No fast exponential deviation inequalities for the progressive mixture rule}

\titlerunning{No fast exponential deviation inequalities in model selection aggregation}

\author{Jean-Yves Audibert}

\authorrunning{J.-Y. Audibert}

\institute{CERTIS - Ecole des Ponts\\
19, rue Alfred Nobel - Cit\'e Descartes\\
77455 Marne-la-Vall\'ee - France\\
\email{audibert@certis.enpc.fr}}

\maketitle

\def\vsp{\vspace{1cm}}
\def\lhs{\text{l.h.s.}}
\def\rhs{\text{r.h.s.}}
\def\wrt{\text{w.r.t.}} 
\def\argmax{\text{argmax}}
\def\argmin{\text{argmin}}

\def\ds1{{\bf 1}} %\mathds{1}}

\def\demi{\frac{1}{2}}
\def\demic{1/2}

\def\E{\mathbb{E}}
\def\N{\mathbb{N}}
\def\P{\mathbb{P}}
\def\R{\mathbb{R}}

\def\begar{$$\begin{array}{lll}}
\def\endar{\end{array}$$}
\def\begarc{$}
\def\endarc{$ }
\def\lbegar{$$\left\{ \begin{array}{lll}}
\def\rendar{\end{array} \right.$$}
\def\rendarp{\end{array} \right..$$}
\def\begarlab{\begin{equation} \begin{array}{lll} \label}
\def\endarlab{\end{array} \end{equation}}
\def\lbegarlab{\begin{equation} \left\{ \begin{array}{lll} \label}
\def\rendarlab{\end{array} \right. \end{equation}}
\def\rendarplab{\end{array} \right.. \end{equation}}

\def\EP{\E_\P}

\def\A{\mathcal{A}}
\def\B{\mathcal{B}}
\def\C{\mathcal{C}}
\def\D{\mathcal{D}}
\def\F{\mathcal{F}}
\def\G{\mathcal{G}}
\def\L{\mathcal{L}}
\def\M{\mathcal{M}}
\def\U{\mathcal{U}}
\def\W{\mathcal{W}}
\def\X{\mathcal{X}}
\def\Y{\mathcal{Y}}
\def\Z{\mathcal{Z}}

\newcommand\und[2]{\underset{#2}{#1}\,}
\newcommand\undc[2]{{#1}_{#2}\;}

\newcommand\fracl[2]{{(#1)}/{#2}}
\newcommand\fracc[2]{{#1}/{#2}}
\newcommand\fracr[2]{{#1}/{(#2)}}
\newcommand\fracb[2]{{(#1)}/{(#2)}}

\def\eqdef{\triangleq}
\def\eps{\epsilon}
\def\logeps{\log(\eps^{-1})}
\def\logeeps{\log(e\eps^{-1})}
\def\lam{\lambda}

\def\calE{\mathcal{E}}
\def\calP{\mathcal{P}}
\def\calR{\mathcal{R}}

\def\pto{\text{o}}
\def\gdo{\text{O}}

\def\hg{\hat{g}}
\def\hh{\hat{h}}
\def\hy{\hat{y}}

\newcommand\integ[1]{\left\lfloor{#1}\rfloor\right.}

\newcommand\refp[1]{\ref{#1} [p.\pageref{#1}]}
\newcommand\eqrefp[1]{\eqref{#1} [p.\pageref{#1}]}

\newcommand\expe[2]{\undc{\E}{#1\sim#2}}
\newcommand\expecv[2]{\E_{#1}}
\newcommand\expecd[2]{\E_{#2}}

\def\hpi{\hat{\pi}}
\def\hrho{\hat{\rho}}

\def\pirho{\hpi({\rho})}
\def\id{\text{I}}
\def\tildY{\tilde{\Y}}
\def\tildg{\tilde{g}}
\def\tildj{J}
\def\tildR{\tilde{R}}

\def\bbb{\text{pim}}
\def\aaa{\text{pm}}
\def\erm{\text{erm}}

\def\pilam{\hat{\pi}_{i}}

\def\tm{\tilde{m}}

\def\ha{h_1}
\def\hb{h_2}
\def\da{d_{\textnormal{I}}}
\def\db{d_{\textnormal{II}}}
\def\tdb{\tilde{\db}}

\def\ya{y_1}
\def\tya{\tilde{\ya}}
\def\yb{y_2}
\def\tyb{\tilde{\yb}}
\def\tg{\tilde{g}}
\def\tga{g_1}
\def\tgb{g_2}
\def\ga{g_1}
\def\gb{g_2}
\def\gz{\tga}
\def\gu{\tgb}

\def\yz{\ya}
\def\yu{\yb}

\def\hri{R(\tg)}
\def\tilam{\tilde{\lam}}

\begin{abstract}
We consider the learning task consisting in predicting as well as the
best function in a finite reference set $\G$ up to the smallest possible additive term.
If $R(g)$ denotes the generalization error of a prediction function $g$,
under reasonable assumptions on the loss function
(typically satisfied by the least square loss when the output is bounded), it is known that 
the progressive mixture rule $\hg$ satisfies 
	\begarlab{eq:1}
	\E R(\hg) \le \undc{\min}{g\in\G} R(g) + C \frac{\log|\G|}{n}, %\log|\G|
	\endarlab
where $n$ denotes the size of the training set, $\E$ denotes the expectation $\wrt$ the
training set distribution and $C$ denotes a positive constant.

This work mainly shows that for any training set size $n$, 
there exist $\eps>0$, a reference set $\G$ and a probability distribution generating the data
such that with probability at least $\eps$ 
	\begar
	R(\hg) \ge \undc{\min}{g\in\G} R(g) + c \sqrt{\frac{\log(|\G|\eps^{-1})}{n}},
	\endar
where c is a positive constant.
In other words, surprisingly, for appropriate reference set $\G$, the deviation convergence rate of the progressive mixture rule is 
only of order $1/\sqrt{n}$ while
its expectation convergence rate is of order $1/n$. The same conclusion
holds for the progressive indirect mixture rule. 
This work also emphasizes on the suboptimality of algorithms based on penalized empirical risk minimization on $\G$.
%the algorithm proposed in see \cite[Section 4]{Aud06a}.
\end{abstract}

\section{Setup and notation}

We assume that we observe $n$ pairs of input-output denoted $Z_1=(X_1,Y_1),\dots,$ $Z_n=(X_n,Y_n)$ 
and that each pair has been independently drawn from the same unknown distribution denoted $P$. 
%These $n$ pairs define the training set and are denoted $Z_1^n$ for short.
The input and output space
are denoted respectively $\X$ and $\Y$, so that $P$ is a probability distribution on the product space 
$\Z \eqdef \X\times\Y$. 
%We assume that $\Y$ is an interval of the real space $\R$.
%To shorten notation, $Z_1^n$ denotes the training set and $\P(dZ_1^n)$ denotes its distribution.
%The target of a learning algorithm is to predict the output $Y$ associated to an input $X$
%for pairs $(X,Y)$ drawn from the distribution $P$. 
The quality of a (prediction) function 
$g:\X\rightarrow\Y$ is measured by the \emph{risk} (or generalization error):
	$$R(g) = {\expe{(X,Y)}{P}} \ell[Y,g(X)],$$
where $\ell[Y,g(X)]$ denotes the loss (possibly infinite) incurred by predicting $g(X)$ when the true output is $Y$.
We work under the following assumptions for the data space 
and the loss function $\ell:\Y\times\Y\rightarrow \R\cup\{+\infty\}.$

\vspace{0.1cm}
\noindent{\bf Main assumptions. }
The input space is assumed to be infinite: $|\X|=+\infty.$
The output space is a non-trivial (i.e. infinite) interval of $\R$ symmetrical $\wrt$ some $a\in\R$: 
for any $y\in\Y$, we have $2a-y\in\Y$.
The loss function is 
\begin{itemize}
\item \emph{uniformly exp-concave:} there exists $\lam>0$ such that
for any $y\in\Y$, the set $\big\{y'\in\R:\ell(y,y')<+\infty \big\}$ is an interval containing $a$ on which
the function $y'\mapsto e^{-\lam \ell(y,y')}$ is concave.
\item \emph{symmetrical:} for any $y_1,y_2\in\Y$,
	\begarc
	\ell(y_1,y_2) = \ell(2a-y_1,2a-y_2),
	\endarc
\item \emph{admissible:} for any $y,y'\in\Y\cap]a;+\infty[$, %and any $y'\in\Y\cap]a;+\infty[$
    $
    \ell(y,2a-y') > \ell(y,y'),
    $
\item \emph{well behaved at center:} for any $y\in\Y\cap]a;+\infty[$, 
the function $\ell_{y}:y' \mapsto \ell(y,y')$ is twice continuously differentiable on a neighborhood of $a$
and $\ell_{y}'(a)<0$.
\end{itemize}

%\begin{remark}
\noindent These assumptions imply that 
\begin{itemize}
\item $\Y$ has necessarily one of the following form: $]-\infty;+\infty[$, $[a-\zeta;a+\zeta]$ or 
$]a-\zeta;a+\zeta[$ for some $\zeta>0$.
\item for any $y\in\Y$, from the exp-concavity assumption, the function $\ell_{y}: y' \mapsto \ell(y,y')$ is convex on the interval on
which it is finite\footnote{Indeed, if $\xi$ denotes the function $e^{-\lam \ell_y}$, from Jensen's inequality,
for any probability distribution, %on $\Y$, 
	$\E \ell_y(Y)=\E\big(-\frac{1}{\lam} \log \xi(Y)\big) 
		\ge -\frac{1}{\lam} \log \E \xi(Y)
		\ge -\frac{1}{\lam} \log \xi(\E Y) = \ell_y(\E Y).$}.
As a consequence, the risk $R$ is also a convex function (on the convex set of prediction functions
for which it is finite).
\end{itemize}
%\end{remark}

\noindent The assumptions were motivated by the fact that they are satisfied in the following settings:
%\vspace{0.1cm}
%\noindent{\bf Examples. }
%The previous assumptions are satisfied in the following settings :
\begin{itemize}
\item least square loss with bounded outputs: $\Y=[y_{\min};y_{\max}]$ and $\ell(y_1,y_2) = (y_1-y_2)^2$. 
Then we have $a=(y_{\min}+y_{\max})/2$ and may take $\lam=1/[2(y_{\max}-y_{\min})^2]$.
\item entropy loss: $\Y=[0;1]$ and
	$\ell(y_1,y_2) = y_1\log\big(\frac{y_1}{y_2}\big)+(1-y_1)\log\big(\frac{1-y_1}{1-y_2}\big)$.
Note that $\ell(0,1) = \ell(1,0) = +\infty$. 
Then we have $a=1/2$ and may take $\lam=1$.
\item exponential (or AdaBoost) loss: $\Y=[-y_{\max};y_{\max}]$ and $\ell(y_1,y_2)=e^{-y_1y_2}$.
Then we have $a=0$ and may take $\lam=e^{-y_{\max}^2}$.
\item logit loss: $\Y=[-y_{\max};y_{\max}]$ and $\ell(y_1,y_2)=\log(1+e^{-y_1y_2})$.
Then we have $a=0$ and may take $\lam=e^{-y_{\max}^2}$.
\end{itemize}

\vspace{0.1cm}
\noindent {\bf Progressive indirect mixture rule.}
Let $\G$ be a finite reference set of prediction functions.
Under the previous assumptions, the only known algorithms satisfying 
\eqref{eq:1}	
are the progressive indirect mixture rules defined below.

For any $i\in\{0,\dots,n\}$, the \emph{cumulative loss} suffered by the prediction function $g$ 
on the first $i$ pairs of input-output is %, denoted $Z_1^i$ for short, is 
	\begar
	\Sigma_i(g) \eqdef \sum_{j=1}^i \ell[Y_j,g(X_j)],
	\endar
where by convention we take $\Sigma_0\equiv 0$.
Let $\pi$ denote 
the uniform distribution on $\G$.
We define the probability distribution $\pilam$ on $\G$ as 
	\begar
	\pilam \propto e^{-\lam \Sigma_i} \cdot \pi
	\endar
equivalently for any $g\in\G$,
	\begarc
	\pilam(g) = \fracr{e^{-\lam \Sigma_i(g)}}{\sum_{g'\in\G} e^{-\lam \Sigma_i(g')}}.
	\endarc
This distribution concentrates on functions having low cumulative loss up to time~$i$.
For any $i\in\{0,\dots,n\}$, let $\hh_i$ be a prediction function such that 
	\begarlab{eq:vov}
	\forall \, (x,y) \in\Z \qquad
		\ell[y,\hh_i(x)] \le -\frac{1}{\lam} \log {\expe{g}{\pilam}} e^{-\lam \ell[y,g(x)]}.
	\endarlab
The \emph{progressive indirect mixture rule}
produces the prediction function 
	\begar
	\hg_{\bbb} = \frac{1}{n+1} \sum_{i=0}^n \hh_i.
	\endar

From the uniform exp-concavity assumption and Jensen's inequality, $\hh_i$ does exist
since one may take $\hh_i= {\expe{g}{\pilam}} g$. This particular choice leads to the 
\emph{progressive mixture rule}, for which the predicted output for any $x\in\X$ is 
	\begar
	\hg_{\aaa}(x) = \sum_{g\in\G} \Big( \frac{1}{n+1} \sum_{i=0}^n 
		\frac{e^{-\lam \Sigma_i(g)}}{\sum_{g'\in\G} e^{-\lam \Sigma_i(g')}} \Big) \, g(x).
	\endar
Consequently, any result that holds for any progressive indirect mixture rule 
in particular holds for the progressive mixture rule.

The idea of a progressive mean of estimators has been introduced by Barron (\cite{Bar87}) in the context
of density estimation with Kullback-Leibler loss. The form $\hg_{\aaa}$ is due to Catoni (\cite{Cat97}). It was also independently proposed in
\cite{Bar99}. The study of this procedure was made in density estimation and least square regression
in \cite{Cat99,Bla99,Yan00,Bun05}. Results for general losses can be found in \cite{Jud06,Aud06a}.
Finally, the progressive indirect mixture rule is inspired by the work 
of Vovk, Haussler, Kivinen and Warmuth \cite{Vov90,Hau98,Vov98} on sequential 
prediction and was studied in the ``batch'' setting in~\cite{Aud06a}.

%Finally, the symbol $\eqdef$ is used to underline when an equality is a definition.
%The symbol $\equiv$ is used when a function is identical to a constant.
%With slight abuse, a symbol denoting a constant function may be used to 
%denote the value of this function.
The symbol $C$ will denote some positive constant whose value may differ
from line to line. %The set of non-negative real numbers is denoted $\Rp=[0;+\infty[$.
The logarithm in base 2 is denoted by $\log_2$ (i.e. $\log_2 t = \log t / \log 2$)
and $\integ{x}$ denotes the largest integer $k$ such that $k \le x$.

\section{Expectation convergence rate}

First let us define the expectation convergence rate of a learning algorithm. 

\begin{definition}
For a given reference set $\G$ of 
prediction functions and a set $\calP$ of probability distributions on $\Z=\X\times\Y$,
a positive sequence $(\Delta_n)_{n\ge 2}$ is said to be an \emph{expectation convergence rate} 
of a learning algorithm producing the prediction function $\hg$ iff 
there exist $C>c>0$ such that
\begin{enumerate}
\item for any distribution $P\in\calP$ and any $n\ge 2$, we have
	\begarlab{eq:exp}
	\E R(\hg) - \inf_{g\in\G} \, R(g) \le C \Delta_n
	\endarlab
\item for large enough $n$, there exists $P\in\calP$ for which
	\begar
	\E R(\hg) - \inf_{g\in\G} \, R(g) \ge c \Delta_n.
	\endar
\end{enumerate}
We say that 
the rate $\Delta_n$ is \emph{optimal}
iff the previous item $2$ is also satisfied for any other algorithm, in other words iff
there is no algorithm having an expectation convergence rate $\tilde{\Delta}_n$
satisfying $\lim_{n\rightarrow +\infty} \tilde{\Delta}_n / {\Delta}_n = 0$.
\end{definition}

The following theorem shows that the expectation convergence rate of any progressive indirect mixture rule is
at least $\fracl{\log |\G|}{n}$ and that for any positive integer $d$, 
there exists a set $\G$ of $d$ prediction functions such that this rate is optimal 
whether we take $\calP$ as the set of all probability distributions on $\Z$
or the set of all probability distributions on $\Z$ for which the output has almost surely two symmetrical
values (e.g. \{-1;+1\}-classication with exponential or logit losses).

\begin{theorem} \label{th:1}
%Let $\hg$ denote the prediction function produced by a progressive indirect mixture rule. We have
Any progressive indirect mixture rule satisfies 
	\begar
	\E R(\hg_{\bbb}) \le \und{\min}{g\in\G} R(g) + \frac{\log|\G|}{\lam(n+1)}.
	\endar
Let $\ya\in\Y-\{a\}$ and $d$ be a positive integer.
There exists a set $\G$ of $d$ prediction functions such that:
for any learning algorithm, 
there exists a probability distribution %on the data space $\X\times[0;1]$ 
generating the data for which
	\begin{itemize}
	\item the output marginal is supported by $2a-\ya$ and $\ya$: $P(Y\in\{ 2a-\ya;\ya \}) = 1$,
	\item 
	\begarc
	\E R(\hg) \ge \und{\min}{g\in\G} R(g) + e^{-1} 
		\kappa \big( 1 \wedge \frac{\lfloor \log_2 |\G| \rfloor}{n+1} \big),
	\endarc
	with $\kappa\eqdef \und{\sup}{y\in\Y} [\ell(\ya,a)-\ell(\ya,y)]>0$.
	\end{itemize}
\end{theorem}

\begin{proof}
See Appendix \ref{sec:proof1}.
\end{proof}

The second part of Theorem \ref{th:1} has the same $(\log |\G|/n)$-rate as the lower bounds obtained in sequential prediction (\cite{Hau98}).
From the link between sequential predictions and our ``batch'' setting with i.i.d. data (see e.g. \cite[Lemma 3]{Aud06a}),
upper bounds for sequential prediction lead to upper bounds for i.i.d. data, and 
lower bounds for i.i.d. data leads to lower bounds for sequential prediction.
The converse of this last assertion is not true, so that the second part of Theorem \ref{th:1}
is not a consequence of the lower bounds of \cite{Hau98}.

The following theorem shows that for appropriate set $\G$:
\begin{itemize}
\item the empirical risk minimizer has a $\sqrt{\fracc{\log|\G|}{n}}$-expectation convergence rate.
\item any empirical risk minimizer and any of its penalized variants are really poor algorithms 
in our learning task since their expectation convergence rate cannot be faster than $\sqrt{\fracc{\log|\G|}{n}}$.
This last point explains the interest we have in progressive mixture rules.
\end{itemize}

\begin{theorem} \label{th:2}
If $B\eqdef \sup_{y,y',y''\in\Y} [ \ell(y,y') - \ell(y,y'') ]<+\infty$, then 
any empirical risk minimizer, which produces a prediction function $\hg_{\erm}$ in $\argmin_{g\in\G} \, \Sigma_n$,
satisfies:
	\begar
	\E R(\hg_\erm) \le \und{\min}{g\in\G} R(g) + B \sqrt{\frac{2\log|\G|}{n}}.
	\endar
Let $\ya,\tya\in\Y\cap]a;+\infty[$ and $d$ be a positive integer.
There exists a set $\G$ of $d$ prediction functions such that:
for any learning algorithm producing a prediction function in $\G$ (e.g. $\hg_\erm$) 
there exists a probability distribution generating the data for which
	\begin{itemize}
	\item the output marginal is supported by $2a-\ya$ and $\ya$: $P(Y\in\{ 2a-\ya;\ya \}) = 1$,
	\item 
	\begarc
	\E R(\hg) \ge \und{\min}{g\in\G} R(g) + \frac{\delta}{8} 
		\Big( \sqrt{\frac{\lfloor \log_2 |\G| \rfloor}{n}} \wedge 2 \Big),
	\endarc
	with $\delta\eqdef \ell(\ya,2a-\tya)-\ell(\ya,\tya)>0$.
	\end{itemize}
\end{theorem}
 
\begin{proof}
See Appendix \ref{sec:proof2}.
\end{proof}

\section{Deviation convergence rate}

The efficiency of an algorithm $\hg$ can be 
summarized by its expected risk $\E \, R(\hg)$, but this does not precise the fluctuations of
$R(\hg)$. In several application fields of learning algorithms, these fluctuations play a key role: 
in finance for instance, the bigger the losses can be, the more money the bank needs to freeze 
in order to alleviate these possible losses. In this case, a ``good'' algorithm is 
an algorithm having not only low expected risk but also small deviations. 

%Let us now define the deviation convergence rate of a learning algorithm. %for a reference set $\G$ of
The deviation convergence rate we define now
is concerned with exponential deviation inequalities (such as Hoeffding's inequality or more generally
such as standard statistical learning inequalities on the supremum of empirical processes).

%\vspace{0.1cm}
%\noindent {\bf Expectation convergence rate and deviation convergence rate.}
\begin{definition}
%Let $\calP$ be a known set of probability distributions on $\Z$ 
%in which we assume that the distribution generating the data is.
%Typical examples for $\calP$ are the set of all probability distributions on $\Z$ 
%and for binary classification with exponential or logit loss
%the set of all probability distributions on $\Z$ for which $Y\in\{-1;+1\}$ almost surely.
%Consider a learning algorithm producing the prediction function $\hg$.
%In the following inequalities, $\E$ denotes the expectation wrt the training set distribution.
Let $0< \gamma \le 1$.
For a given reference set $\G$ of 
prediction functions and a set $\calP$ of probability distributions on $\Z=\X\times\Y$,
a positive sequence $(\Delta'_n)_{n\in\N}$ is said to be a \emph{deviation convergence rate
of order $\gamma$} 
of a learning algorithm iff there exist $C>c>0$ such that
\begin{enumerate}
\item for any distribution $P\in\calP$, integer ${n\ge2}$, and $\eps>0$, with probability at least $1-\eps$ 
$\wrt$ the training set distribution, we have
	\begarlab{eq:dev}
	R(\hg) - \inf_{g\in\G} \, R(g) \le C \big[\log^{\gamma}(e\eps^{-1})\big] \Delta'_n,
	\endarlab
\item for large enough $n$, there exist $\eps>0$ and a distribution $P\in\calP$ such that
with probability at least $\eps$ $\wrt$ the training set distribution, we have
	\begar
	R(\hg) - \inf_{g\in\G} \, R(g) \ge c \big[\log^{\gamma}(e \eps^{-1})\big] \Delta'_n.
	\endar
\end{enumerate}
\end{definition}

The following lemma shows that the expectation convergence rate 
of a learning algorithm is at least of order of its deviation convergence rate.
The expectation convergence rate can also be strictly faster as the comparison between
Theorems \ref{th:1} and \ref{th:main1} shows.

\begin{lemma}
Let $\hg$ satisfy: for any $\eps>0$, with probability at least $1-\eps$, \eqref{eq:dev} holds.
Then we have 
	\begar
	\E R(\hg) - \inf_{g\in\G} \, R(g) \le 2^{\gamma} C \Delta'_n.
	\endar
\end{lemma}

\begin{proof}
It suffices to integrate the deviations. Let $R^*= \inf_{g\in\G} \, R(g).$
By Jensen's inequality, we have
	\begar
	\Big(\frac{\E R(\hg) - R^*}{C \Delta'_n}\Big)^{1/\gamma} - 1\\
	\qquad
	\le \E \Big(\frac{R(\hg) - R^*}{C \Delta'_n}\Big)^{1/\gamma} - 1\\
	\qquad
		\le \E \Big\{ \Big[\big(\frac{R(\hg) - R^*}{C \Delta'_n}\big)^{1/\gamma} - 1 \Big] \vee 0 \Big\}\\
	\qquad
		= \int_{0}^{+\infty} \P\big\{ 
		\big(\frac{R(\hg) - R^*}{C \Delta'_n}\big)^{1/\gamma} - 1 > u \big\} du\\
	\qquad
		= \int_{0}^{1} \P\big\{ R(\hg) - R^* > C \Delta_n'\log^{\gamma}(e\eps^{-1}) \big\} \frac{d\eps}{\eps}
		\qquad \text{[setting $u=\logeps$]}\\
	\qquad
		\le 1.
	\endar
\end{proof}

The following theorem shows that the deviation convergence rate of order $1/2$ of
any progressive indirect mixture rule is at least %$\sqrt{\fracl{\log |\G|}{n}}$ 
$1/\sqrt{n}$ and that
there exists $\G$ such that the deviation convergence rate of order $1/2$ of any progressive indirect mixture rule is 
$1/\sqrt{n}$ %(up to possibly a logarithmic factor)
whether we take $\calP$ as the set of all probability distributions on $\Z$
or the set of all probability distributions on $\Z$ for which the output has almost surely two symmetrical
values (e.g. \{-1;+1\}-classication with exponential or logit losses).

\begin{theorem} \label{th:main1}
%Let $\hg$ denote the prediction function produced by a progressive indirect mixture rule. We have
If $B\eqdef \sup_{y,y',y''\in\Y} [ \ell(y,y') - \ell(y,y'') ]<+\infty$, then 
any progressive indirect mixture rule satisfies: %for $n$ large enough (see \eqrefp{eq:nlargeenough}), 
%for any $n$, 
for any $\eps>0$, with probability at least $1-\eps$ $\wrt$ the training set distribution, we have
	\begar
	R(\hg_{\bbb}) \le \und{\min}{g\in\G} R(g) 
	 + B\sqrt{\frac{2\logeps}{n+1}} + \frac{\log |\G|}{\lam(n+1)} 
	\endar
Let $\ya$ and $\tya$ in $\Y\cap ]a;+\infty[$ such that 
$\ell_{\ya}$ is twice continuously differentiable on $[a;\tya]$ and
$\ell_{\ya}'(\tya)\le 0$ and $\ell_{\ya}''(\tya) > 0$.
Consider the prediction functions $\tga\equiv\tya$ and $\tgb\equiv 2a-\tya$.
%For any prior distribution $\pi$ and positive parameter $\lam$,
%Whether $g$ denotes $\hg_{\aaa}$ or $\hg_{\bbb}$,
%For any probability distribution $\pi$ and positive parameter $\lam$, there exists $n_0$,
For any training set size $n$ large enough, there exist $\eps>0$ and a distribution
generating the data such that
\begin{itemize}
\item the output marginal is supported by $\ya$ and $2a-\ya$
\item with probability larger than $\eps$, we have
	\begar
	R(\hg_{\bbb})- \und{\min}{g\in\{\tga,\tgb\}} R(g) & \ge & c \sqrt{\frac{\log(e\eps^{-1})}{n}}
	\endar
where $c$ is a positive constant depending only on the loss function, the symmetry parameter $a$ and the output values $\ya$ and $\tya.$ 
\end{itemize}
\end{theorem}

\begin{proof}
See Section \ref{sec:proofthmain1}.
\end{proof}

This result is quite surprising since it gives an example of an algorithm which is optimal
in terms of expectation convergence rate and for which the deviation convergence rate is 
(significantly) worse that the expectation convergence rate. 
%One of the result of \cite{Aud07b} is to prove that the deviation convergence rate of the
%progressive indirect mixture rule is not optimal to the extent that for any $\eps>0$,
%there exists a learning algorithm producing a prediction function $\hg$ such that TMP.

\section{Proof of Theorem \ref{th:main1}} \label{sec:proofthmain1}

\subsection{Proof of the upper bound} 

We would like to thank an anonymous reviewer for suggesting the following proof, which leads to better constants
than the original one based on PAC-Bayesian inequalities.

Let $Z_{n+1}=(X_{n+1},Y_{n+1})$ be an input-output pair
independent from the training set $Z_1,\dots,Z_n$ and with the same distribution $P$.
From the convexity of $y' \mapsto \ell(y,y')$, we have 
	\begarlab{eq:r1}
	R(\hg_{\bbb}) \le \frac{1}{n+1} \sum_{i=0}^n R(\hh_i).
	\endarlab
Now from \cite[Theorem 1]{Zha05} (see also \cite[Proposition 1]{Ces04}), 
for any $\eps>0$, with probability at least $1-\eps$, we have
	\begarlab{eq:r2}
	\frac{1}{n+1} \sum_{i=0}^n R(\hh_i) \le \frac{1}{n+1} \sum_{i=0}^n \ell\big(Y_{i+1},\hh(X_{i+1})\big)
		+ B\sqrt{\frac{\logeps}{2(n+1)}}
	\endarlab
Using \cite[Theorem 3.8]{Hau98} and the exp-concavity assumption, we have
	\begarlab{eq:r3}
	\sum_{i=0}^n \ell\big(Y_{i+1},\hh(X_{i+1})\big) \le 
		\und{\min}{g\in\G} \sum_{i=0}^n \ell\big(Y_{i+1},g(X_{i+1})\big) + \frac{\log |\G|}{\lam}		
	\endarlab
Let $\tg \in \argmin_{\G} \, R$. By Hoeffding's inequality, with probability at least $1-\eps$, we have
	\begarlab{eq:r4}
	\frac{1}{n+1} \sum_{i=0}^n \ell\big(Y_{i+1},\tg(X_{i+1})\big) \le R(\tg)+ B\sqrt{\frac{\logeps}{2(n+1)}}
	\endarlab
Merging \eqref{eq:r1}, \eqref{eq:r2}, \eqref{eq:r3} and \eqref{eq:r4}, 
%for any $\eps>0$, 
with probability at least $1-2\eps$, we get
	\begar
	R(\hg_{\bbb}) & \le & \frac{1}{n+1} 
		\sum_{i=0}^n \ell\big(Y_{i+1},\tg(X_{i+1})\big) + \frac{\log |\G|}{\lam(n+1)} 
		+ B\sqrt{\frac{\logeps}{2(n+1)}}\\
	& \le & R(\tg) + B\sqrt{\frac{2\logeps}{n+1}} + \frac{\log |\G|}{\lam(n+1)}.
	\endar

\subsection{Proof of the lower bound} 

We cannot use standard tools like Assouad's argument (see e.g. \cite[Theorem 14.6]{Dev96})
because if it were possible, it would mean that the lower bound would hold for any algorithm 
and this is (non trivially) false. %see TMP

To prove that any progressive indirect mixture rule have no fast exponential deviation inequalities, we will 
show that on some event with not too small probability, for most of the $i$ in $\{0,\dots,n\}$,
$\pi_{-\lam \Sigma_i}$ concentrates on the wrong function. %, that is $\gb$.

%Now let us give the technical details of the proof. 
The proof is organized as follows. First we
define the probability distribution for which we will prove that the progressive indirect mixture rules cannot have
fast deviation convergence rates. Then we define the event on which the progressive indirect 
mixture rules do not perform well.
We lower bound the probability of this excursion event.
Finally we conclude by lower bounding $R(\hg_{\bbb})$ on the excursion event.

Before starting the proof, note that from the ``well behaved at center'' and exp-concavity assumptions,  
for any $y\in\Y\cap]a;+\infty[$, on a neighborhood of $a$, we have: $\ell_y''\ge\lam(\ell_y')^2$ and since $\ell_y'(a)<0$,
$\ya$ and $\tya$ exist.

\subsubsection{Probability distribution generating the data and first consequences.}

Let $\gamma\in]0;1]$ be a parameter to be tuned later.
We consider a distribution generating the data such that the output distribution satisfies for any $x\in\X$
	\begar
	P(Y=\ya|X=x)=\fracl{1+\gamma}{2}=1-P(Y=\yb|X=x),
	\endar
where $\yb=2a-\ya.$
Let $\tyb=2a-\tya.$
From the symmetry and admissibility assumptions, we have
	\begarc
	\ell(\yb,\tyb) = \ell(\ya,\tya) < \ell(\ya,\tyb) = \ell(\yb,\tya).
	\endarc
Introduce
	\begarlab{eq:defdelta}
	\delta \eqdef \ell(\ya,\tyb)-\ell(\ya,\tya)>0.
	\endarlab 
We have
	\begarlab{eq:riskgap}
	R(\gb)-R(\ga) = \frac{1+\gamma}{2}[\ell(\ya,\tyb)-\ell(\ya,\tya)]
		+ \frac{1-\gamma}{2}[\ell(\yb,\tyb)-\ell(\yb,\tya)]
		= \gamma \delta. 
	\endarlab
Therefore $\ga$ is the best prediction function in $\{\ga,\gb\}$ for the distribution we have chosen.
Introduce $W_j\eqdef \ds1_{Y_j=\ya}-\ds1_{Y_j=\yb}$
and $S_i \eqdef \sum_{j=1}^i W_j$. 
For any $i\in\{1,\dots,n\}$, we have 
	\begar
	\Sigma_i(\gb)-\Sigma_i(\ga) = \sum_{j=1}^i [\ell(Y_j,\tyb)-\ell(Y_j,\tya)]
	    = \sum_{j=1}^i W_j\delta = \delta \, S_i
	\endar
The weight given by the Gibbs distribution $\pi_{-\lam \Sigma_i}$ to 
the function $\ga$ is
	\begarlab{eq:lbpisigmai}
	\pi_{-\lam \Sigma_i}(\ga) 
		= \frac{e^{-\lam \Sigma_i(\ga)}}{e^{-\lam \Sigma_i(\ga)}+e^{-\lam \Sigma_i(\gb)}}
		= \frac{1}{1+e^{\lam [\Sigma_i(\ga)-\Sigma_i(\gb)]}}
		= \frac{1}{1+e^{-\lam \delta S_i}}.
	\endarlab	
	
\subsubsection{An excursion event on which the progressive indirect mixture rules will not perform well.}
\eqref{eq:lbpisigmai} leads us to consider the event: 
	\begar
	E_\tau = \big\{ \forall i\in\{\tau,\dots,n\}, \; S_i \le - \tau \big\},
	\endar
with $\tau$ the smallest integer larger than $(\log n)/(\lam \delta)$
such that $n-\tau$ is even. (We could have just as well chosen $n-\tau$ odd; see \eqref{eq:low1a} below.)
We have
	\begarlab{eq:deftau}
	\frac{\log n}{\lam\delta}\le \tau \le \frac{\log n}{\lam\delta}+2.
	\endarlab

The event $E_\tau$ can be seen as an excursion event of the random walk defined through 
the random variables $W_j=\ds1_{Y_j=\ya}-\ds1_{Y_j=\yb}$, $j\in\{1,\dots,n\}$, which are equal to 
$+1$ with probability $(1+\gamma)/2$ and $-1$ with probability $(1-\gamma)/2$.

From \eqref{eq:lbpisigmai}, on the event $E_\tau$, for any $i\in\{\tau,\dots,n\}$, we have 
	\begarlab{eq:badconc}
	\pi_{-\lam \Sigma_i}(\ga) \le \frac{1}{n+1}.
	\endarlab
This means that %on this event for any $i\ge \tau$, 
$\pi_{-\lam \Sigma_i}$ concentrates on the wrong function,
i.e. the function $g_2$ having larger risk (see \eqref{eq:riskgap}).
	
\subsubsection{Lower bound of the probability of the excursion event.} 

This requires to look at the probability
that a slightly shifted random walk in the integer space has a very long excursion
above a certain threshold.
To lower bound this probability, we will first look at the non-shifted random walk.
Then we will see that for small enough shift parameter, probabilities of shifted random walk events are close to 
the ones associated to the non-shifted random walk.

Let $N$ be a positive integer.
Let $\sigma_1,\dots,\sigma_N$ be $N$ independent Rademacher variables: $\P(\sigma_i=+1)=\P(\sigma_i=-1)=1/2$.
Let $s_i \eqdef \sum_{j=1}^i \sigma_i$
be the sum of the first $i$ Rademacher variables.
We start with the following lemma for sums of Rademacher variables.
\begin{lemma} \label{lem:l1}
%Let $m>0$ and $\tau$ be an integer such that $N-\tau$ is even. We have
Let $m$ and $t$ be positive integers. We have
	\begarlab{eq:l1}
	\P\big( \und{\max}{1\le k \le N} s_k \ge t; s_N \neq t ; \big| s_N - t \big| \le m \big)
	 	= 2\P\big( t < s_N \le t + m \big)
	\endarlab
\end{lemma}

\begin{proof}[of Lemma \ref{lem:l1}]
The result comes from the well known mirror trick used to compute the law of 
$\big({\sup}_{s\le t} W_s, W_t \big)$ where $W$ denotes a Brownian motion.
Consider a sequence $\sigma_1,\dots,\sigma_N$ which belongs to the event $\calE$ of the $\lhs$ probability.
Let $\tildj$ be the first integer $j$ such that $s_j=t$.
Since 
\begin{itemize}
\item the sequences $\sigma_1,\dots,\sigma_N$ and 
$\sigma_1,\dots,\sigma_{\tildj},-\sigma_{\tildj+1},\dots,-\sigma_N$ 
have the same probabilities,
\item both sequences belong to $\calE$ and are different since $J<N$,
\item exactly one of the sequences satisfy $s_N > t$,
\end{itemize}
we have
	\begar
	\P\big( \und{\max}{1\le k \le N} s_k \ge t; s_N \neq t ; \big| s_N - t \big| \le m \big)
	 	= 2 \P\big( s_N > t; \big| s_N - t \big| \le m  \big),
	\endar
which is the desired result.
\end{proof}

Let $\sigma'_1,\dots,\sigma'_N$ be $N$ independent shifted Rademacher variables to the extent that
	$\P(\sigma'_i=+1)=(1+\gamma)/2=1-\P(\sigma'_i=-1)$.
These random variables satisfy the following key lemma

\begin{lemma} \label{lem:l2}
For any set 
	$
	A \subset \big\{ (\eps_1,\dots,\eps_N)\in\{-1,1\}^n: \big| \sum_{i=1}^N \eps_i \big| \le M \big\}
	$
where $M$ is a positive integer,
we have
	\begarlab{eq:l2}
	\P\big\{ (\sigma'_1,\dots,\sigma'_N)\in A \big\} \ge \Big( \frac{1-\gamma}{1+\gamma} \Big)^{M/2}
		\big(1-\gamma^2\big)^{N/2} \P\big\{ (\sigma_1,\dots,\sigma_N)\in A \big\}
	\endarlab
\end{lemma}
\begin{proof}[of Lemma \ref{lem:l2}]
Let $s$ be an integer such that $N-s$ is even and $|s|\le M$
Consider a sequence $\eps_1,\dots,\eps_N$ such that $\sum_{i=1}^N \eps_i = s$.
Then the numbers of $-1$ and $+1$ in the sequence are respectively $(N-s)/2$ and $(N+s)/2$.
Consequently, we have
	\begar
	\frac{\P[ (\sigma'_1,\dots,\sigma'_N)=(\eps_1,\dots,\eps_N) ]}
		{\P[ (\sigma_1,\dots,\sigma_N)=(\eps_1,\dots,\eps_N) ]}
		= (1+\gamma)^{(N-s)/2}(1-\gamma)^{(N+s)/2},
	\endar
hence
	\begar
	\P\{ (\sigma'_1,\dots,\sigma'_N)=(\eps_1,\dots,\eps_N) \big\}\\
	\qquad\qquad\qquad
		\ge (1-\gamma^2)^{N/2} \big( \frac{1-\gamma}{1+\gamma} \big)^{M/2} 
		\P\big\{ (\sigma_1,\dots,\sigma_N)=(\eps_1,\dots,\eps_N) \big\}.
	\endar
By summing over the sequences $\eps_1,\dots,\eps_N$ in $A$, we obtain the desired result.
\end{proof}

We may now lower bound the probability of the excursion event $E_\tau$.
Let $M$ be an integer larger than $\tau$.
We still use $W_j\eqdef \ds1_{Y_j=\ya}-\ds1_{Y_j=\yb}$ for $j\in\{1,\dots,n\}$.
By using Lemma \ref{lem:l2} with $N=n-2\tau$, we obtain
	\begarlab{eq:low3}
	\P(E_\tau) & \ge & \P\big( W_1=-1,\dots,W_{2\tau}=-1; \;\forall \, 2 \tau < i \le n,\;
		\sum_{j=2\tau+1}^i W_j \le \tau \big)\\
	& = & \big(\frac{1-\gamma}{2}\big)^{2\tau}
		\P\big(\forall \, i> 2 \tau \quad \sum_{j=2\tau+1}^i W_j \le \tau \big)\\
	& = & \big(\frac{1-\gamma}{2}\big)^{2\tau}
		\P\big(\forall \, i\in\{1,\dots,N\} \quad \sum_{j=1}^i \sigma'_j \le \tau \big)\\
	& \ge & \big(\frac{1-\gamma}{2}\big)^{2\tau}
		\P\big(\big| \sum_{i=1}^N \sigma'_i \big| < M;\forall \, i\in\{1,\dots,N\} \quad \sum_{j=1}^i \sigma'_j \le \tau \big)\\
	& \ge & \big(\frac{1-\gamma}{2}\big)^{2\tau} \big(\frac{1-\gamma}{1+\gamma}\big)^{M/2} 
		\big(1-\gamma^2\big)^{\frac{N}{2}}
		\P\big(|s_N| \le M ; \forall \, i\in\{1,\dots,N\} \quad s_i \le \tau \big)\\
	\endarlab
By using Lemma \ref{lem:l1}, since $\tau\le M$, the $\rhs$ probability can be lower bounded: 
	{\renewcommand{\arraystretch}{1.5}
	\begar %lab{eq:low2}
	\P\big(|s_N| \le M ; \und{\max}{1\le i \le N} s_i \le \tau \big)\\
	\qquad\qquad
		= \P\Big\{ \und{\max}{1\le i\le N} \, s_i \le \tau ; s_N \ge -M\Big\}\\
	\qquad\qquad
		\ge \P\Big\{ \und{\max}{1\le i\le N} \, s_i < \tau ; |s_N-\tau| \le M+\tau ; s_N \neq \tau \Big\}\\
	\qquad\qquad
		= \P\Big\{ |s_N-\tau| \le M+\tau ; s_N \neq \tau \Big\}\\
	\qquad\qquad\qquad\qquad
		- \P\Big\{ \und{\max}{1\le i\le N} \, s_i \ge \tau ; |s_N-\tau| \le M+\tau ; s_N \neq \tau \Big\}\\ 
	\qquad\qquad
		= \P\big\{ |s_N-\tau| \le M+\tau ; s_N \neq \tau \big\}
		- 2\P\big\{ \tau < s_N \le M+2\tau \big\}\\ 
	\qquad\qquad
		= \P\big\{ -M \le s_N < \tau \big\} - \P\big\{ \tau < s_N \le M+2\tau \big\}\\
	\qquad\qquad
		= \P\big\{ -\tau < s_N \le M \big\} - \P\big\{ \tau < s_N \le M+2\tau \big\}\\
	\qquad\qquad
		= \P\big\{ -\tau < s_N \le \tau \big\} - \P\big\{ M < s_N \le M+2\tau \big\}\\
	\endar\par}
%Since $n-\tau$ is even, we have $N-\tau=n-3\tau$. Similarly, 
Let us consider only the integer $M>\tau$ 
such that $n-M$ is even, or equivalently $N-M$ is even. Since $N-\tau=n-3\tau$ is also even, we have
	\begarlab{eq:low1a}
	\P\big(|s_N| \le M ; \und{\max}{1\le i \le N} s_i \le \tau \big) \\
	\qquad\qquad\qquad	
		\ge \sum_{k=0}^{\tau-1} \P( s_N =2-\tau + 2k )
		- \sum_{k=1}^{\tau} \P( s_N =M + 2k )\\
	\qquad\qquad\qquad
		\ge \tau [ \P( s_N = \tau ) - \P( s_N = M ) ],
	\endarlab
where the last inequality comes from properties of the binomial coefficients.
	
Combining \eqref{eq:low3} and \eqref{eq:low1a}, we obtain
	\begarlab{eq:low4}
	\P(E_\tau) & \ge & \tau \big(\frac{1-\gamma}{2}\big)^{2\tau} \big(\frac{1-\gamma}{1+\gamma}\big)^{M/2} 
		\big(1-\gamma^2\big)^{\frac{N}{2}} [ \P( s_N = \tau ) - \P( s_N = M ) ]
	\endarlab
where we recall that $\tau$ have the order of 
$\log n$, $N=n-2\tau$ has the order of $n$ and that
 $\gamma > 0$ and $M\ge \tau$ have to be appropriately chosen.

To control the probabilities of the $\rhs$, we use Stirling's formula
	\begarlab{eq:stir1}
	n^n e^{-n} \sqrt{2\pi n} \, e^{1/(12n+1)} < n!  < n^n e^{-n} \sqrt{2\pi n} \, e^{1/(12n)},
	\endarlab
and get for any $s\in[0;N]$ such that $N-s$ even,
	\begarlab{eq:stir2}
	\P(s_N=s) & = & \big(\frac{1}{2}\big)^{N} \binom{N}{\frac{N+s}{2}}\\
	& \ge & \big(\frac{1}{2}\big)^{N}  \frac{ (\frac{N}{e})^N \sqrt{2\pi N} e^{\frac{1}{12N+1}} }
		{ (\frac{N+s}{2e})^{\frac{N+s}{2}} (\frac{N-s}{2e})^{\frac{N-s}{2}} \sqrt{\pi(N+s)}\sqrt{\pi(N-s)}
		e^{\frac{1}{6(N+s)}}e^{\frac{1}{6(N-s)}} }	\\
	& = & \frac{ 1 }{(1+\frac{s}{N})^{\frac{N+s}{2}} (1-\frac{s}{N})^{\frac{N-s}{2}}}
		\sqrt{\frac{2N}{\pi(N^2-s^2)}} e^{\frac{1}{12N+1}
		-\frac{1}{6(N+s)}-\frac{1}{6(N-s)}}	\\
	& \ge & \sqrt{\frac{2}{\pi N}} \Big( 1-\frac{s^2}{N^2} \Big)^{-\frac{N}{2}}
		\Big( \frac{ 1-\frac{s}{N} }{ 1+\frac{s}{N} }\Big)^{\frac{s}{2}}
		e^{-\frac{1}{6(N+s)}-\frac{1}{6(N-s)}}	\\
	\endarlab
and similarly 
	\begarlab{eq:stir3}
	\P(s_N=s) \le \sqrt{\frac{2}{\pi N}} \Big( 1-\frac{s^2}{N^2} \Big)^{-\frac{N}{2}}
		\Big( \frac{ 1-\frac{s}{N} }{ 1+\frac{s}{N} }\Big)^{\frac{s}{2}}
		e^{\frac{1}{12N+1}}
	\endarlab

%For the bracketed term of \eqref{eq:low4} not to be small, 
These computations and \eqref{eq:low4} leads us to take $M$ as the smallest integer larger than $\sqrt{n}$
such that $n-M$ is even. Indeed, from \eqref{eq:deftau}, \eqref{eq:stir2} and \eqref{eq:stir3}, we obtain
	\begarc
	\lim_{n\rightarrow +\infty} \sqrt{n} [\P( s_N = \tau ) - \P( s_N = M )] = c,
	\endarc
where $c=\sqrt{\fracc{2}{\pi}}\big(1-{e}^{-1/2}\big)>0$. Therefore for $n$ large enough we have
	\begarlab{eq:low5}
	\P(E_\tau) & \ge & \frac{c \tau}{2\sqrt{n}} \big(\frac{1-\gamma}{2}\big)^{2\tau} \big(\frac{1-\gamma}{1+\gamma}\big)^{M/2} 
		\big(1-\gamma^2\big)^{\frac{N}{2}} 
	\endarlab
The last two terms of the $\rhs$ of \eqref{eq:low5}
leads us to take $\gamma$ of order $1/\sqrt{{n}}$
up to possibly a logarithmic term.
We obtain the following lower bound on 
the excursion probability
\begin{lemma} \label{lem:3}
If $\gamma = \sqrt{C_0 \fracl{\log n}{n}}$ with $C_0$ a positive constant, then for any large enough $n$,
	\begar
	\P(E_\tau) \ge \frac{1}{n^{C_0}}.
	\endar
\end{lemma}

\subsubsection{Behavior of the progressive indirect mixture rule on the excursion event.}

From now on, we work on the event $E_\tau$. 
%Let $\hg_{\bbb}$ be the prediction function produced by Algorithm \bbb. 
We have $\hg_{\bbb} = (\sum_{i=0}^n \hh_i)/(n+1)$. We still use 
	$\delta \eqdef \ell(\ya,\tyb) - \ell(\ya,\tya) = \ell(\yb,\tya) - \ell(\yb,\tyb).$ 
On the event $E_\tau$, for any $x\in\X$ and any $i\in\{\tau,\dots,n\}$, by definition of $\hh_i$, we have
	\begar
	\ell[\yb,\hh_i(x)] - \ell(\yb,\tyb) & \le & -\frac{1}{\lam} \log {\expe{g}{\pi_{-\lam \Sigma_i}}} 
		e^{-\lam \{\ell[\yb,g(x)]- \ell(\yb,\tyb)\}}\\
	& = & -\frac{1}{\lam} \log \big\{ \pi_{-\lam \Sigma_i}(\ga) e^{-\lam \delta}+\pi_{-\lam \Sigma_i}(\gb) \big\}\\
	& = & -\frac{1}{\lam} \log \big\{ e^{-\lam \delta}+(1-e^{-\lam \delta})\pi_{-\lam \Sigma_i}(\gb) \big\}\\
	& \le & -\frac{1}{\lam} \log \big\{ 1-(1-e^{-\lam \delta})\frac{1}{n+1} \big\}\\
	\endar
In particular, for any $n$ large enough, we have 
	$
	\ell[\yb,\hh_i(x)] - \ell(\yb,\tyb) \le Cn^{-1},
	$
with $C>0$ \emph{independent from $\gamma$}.
From the convexity of the function $y\mapsto \ell(\yb,y)$ and 
by Jensen's inequality, we obtain
	\begarlab{eq:gbprop}
	\ell[\yb,\hg_{\bbb}(x)] - \ell(\yb,\tyb) 
		& = & \ell[\yb,\frac{1}{n+1} \sum_{i=0}^n \hh_i(x)] - \ell(\yb,\tyb) \\
	& \le & \frac{1}{n+1} \sum_{i=0}^n \ell[\yb,\hh_i(x)] - \ell(\yb,\tyb) \\
	& \le & \frac{\tau \delta}{n+1} + C n^{-1}\\
	& < & C_1 \frac{\log n}{n}
	\endarlab
for some constant $C_1>0$ \emph{independent from $\gamma$}.
%To prove that the previous inequality implies that $\hg_{\bbb}(x)$ is smaller than $\tyb$ up to 
%some small additive term, we use 
%\begin{lemma}
Let us now prove that for $n$ large enough, we have 
	\begarlab{eq:bbbxprop}
	\tyb \le \hg_{\bbb}(x) \le \tyb+C\sqrt{\frac{\log n}{n}} \le \tya,
	\endarlab
with $C>0$ \emph{independent from $\gamma$}.

%\end{lemma}

\def\tdy{t}
\begin{proof}
For any $y\in\Y$, let $\tdy=2a-y$. 
We have
	$\ell(\yb,y) - \ell(\yb,\tyb) = \ell_{\ya}(\tdy) - \ell_{\ya}(\tya)$.
Since $\ell_{\ya}'(\tya)\le0$, $\ell_{\ya}''(\tya)>0$, $\ell_{\ya}'' \ge \lam (\ell_{\ya}')^2$
and $\ell_{\ya}''$ is continuous on $[a;\tya]$,
there exists $m>0$ such that $\ell_{\ya}'' > m$ on $[a ; \tya].$
For any $\tyb < y \le a$, from Taylor's expansion, we have
	\begarlab{eq:yprop}
	\ell(\yb,y) - \ell(\yb,\tyb) 
	& > & (t-\tya) \ell_{\ya}'(\tya) + \frac{(t-\tya)^2}{2} m \\
	& \ge & \frac{(t-\tya)^2}{2} m \\
	& = & \frac{(y-\tyb)^2}{2} m \\
	\endarlab
Let $y_0 \eqdef \tyb + \sqrt{\frac{2C_1 \log n}{m n}}$ where $C_1$ is the constant appearing in \eqref{eq:gbprop}.
For $n$ large enough, we have $y_0\le a $ and  we may apply \eqref{eq:yprop} to $y=y_0$.
We get 
	\begarlab{eq:y0prop}
	\ell(\yb,y_0) - \ell(\yb,\tyb) > C_1 \frac{\log n}{n}.
	\endarlab
Since $\ell_{\ya}$ is convex, $\ell_{\ya}'(\tya)\le 0$ and $\ell_{\ya}''(\tya)>0$, the function $\ell_{\ya}$ decreases on $]-\infty;\tya] \cap \Y$.
By symmetry, the function $y\mapsto \ell(\yb,y)$ is non-decreasing on $[\tyb;+\infty[ \cap~\Y$.
From \eqref{eq:gbprop} and \eqref{eq:y0prop}, we get $\hg_{\bbb}(x) \notin [y_0;+\infty[,$
which ends the proof of the upper bound of $\hg_{\bbb}(x)$.

For the lower bound, for any $x\in\X$, by definition of $\hh_i$, we have
	\begar
	\ell[\ya,\hh_i(x)] - \ell(\ya,\tya) & \le & -\frac{1}{\lam} \log {\expe{g}{\pi_{-\lam \Sigma_i}}} 
		e^{-\lam \{\ell[\ya,g(x)]- \ell(\ya,\tya)\}}\\
	& = & -\frac{1}{\lam} \log \big\{ \pi_{-\lam \Sigma_i}(\ga) +\pi_{-\lam \Sigma_i}(\gb) e^{-\lam \delta} \big\}\\
	& \le & \delta.
	\endar
By Jensen's inequality, we obtain
	\begar
	\ell_{\ya}[\hg_{\bbb}(x)] - \ell_{\ya}(\tya) 
		& = & \ell[\ya,\frac{1}{n+1} \sum_{i=0}^n \hh_i(x)] - \ell(\ya,\tya)\\
	& \le & \frac{1}{n+1} \sum_{i=0}^n \ell[\ya,\hh_i(x)] - \ell(\ya,\tya)\\
	& \le & \delta\\
	& = & \ell_{\ya}(\tyb ) - \ell_{\ya}( \tya ).
	\endar
Since the function $\ell_{\ya}$ decreases on $]-\infty;\tyb] \cap \Y$, 
we get that $\hg_{\bbb}(x) \ge \tyb$, which ends the proof of \eqref{eq:bbbxprop}.
%\begin{flushright}
%$\blacksquare$
%\end{flushright}
\end{proof}

From \eqref{eq:bbbxprop}, we obtain
	\begarlab{eq:argab}
	R(\hg_{\bbb}) - R( \ga ) & = & \frac{1+\gamma}{2} \big[ \ell(\ya,\hg_{\bbb}) - \ell(\ya,\tya) \big]
		+ \frac{1-\gamma}{2} \big[ \ell(\yb,\hg_{\bbb}) - \ell(\yb,\tya) \big]\\
	& = & \frac{1+\gamma}{2} \big[ \ell_{\ya}(\hg_{\bbb}) - \ell_{\ya}(\tya) \big]
		+ \frac{1-\gamma}{2} \big[ \ell_{\ya}(2a-\hg_{\bbb}) - \ell_{\ya}(\tyb) \big]\\
	& = & \frac{1+\gamma}{2} \big[ \delta + \ell_{\ya}(\hg_{\bbb})-\ell_{\ya}(\tyb) \big]\\
	& & \qquad
		+ \frac{1-\gamma}{2} \big[ -\delta+ \ell_{\ya}(2a-\hg_{\bbb}) - \ell_{\ya}(\tya) \big]\\
	& \ge & \gamma \delta - %\frac{1+\gamma}{2} 
	        (\hg_{\bbb}-\tyb) | \ell_{\ya}'(\tyb) |\\
		%- \frac{1-\gamma}{2} (\hg_{\bbb}-\tyb) | \ell_{\ya}'(\tya) |\\
	& \ge & \gamma \delta - C_2 \sqrt{\frac{\log n}{n}},
	\endarlab
with $C_2$ \emph{independent from $\gamma$}.
We may take $\gamma = \frac{2 C_2}{\delta} \sqrt{\fracl{\log n}{n}}$ and obtain:
for $n$ large enough, on the event $E_\tau$, we have
	\begarc
	R(\hg_{\bbb}) - R( \ga ) 
		\ge C \sqrt{\fracc{\log n}{n}}.
	\endarc
From Lemma \ref{lem:3}, this inequality holds with probability at least 
$1/n^{C_4}$ for some $C_4>0$.
To conclude, for any $n$ large enough, there exists $\eps>0$ s.t. with probability at least 
$\eps$,
	\begar
	R(\hg_{\bbb}) - R( \ga ) 
		\ge c \sqrt{\frac{\log(e\eps^{-1})}{n}}.
	\endar
where $c$ is a positive constant depending only on the loss function, the symmetry parameter $a$ and the output values $\ya$ and $\tya.$ 

\begin{remark}
Had we consider the progressive mixture rule, this last part of the proof
would have been much simpler. Indeed, for $n$ large enough,
on the event $E_\tau$, from \eqref{eq:badconc}, we have 
	\begar %lab{eq:defp}
	p\eqdef \frac{1}{n+1} \sum_{i=0}^n \pi_{-\lam \Sigma_i}(\ga)
		\le \frac{\tau}{n+1}+ \und{\sup}{\tau \le i \le n} \pi_{-\lam \Sigma_i}(\ga)
		\le C \frac{\log n}{n}
	\endar %lab
and 
	\begarc
	\hg_{\aaa} = \frac{1}{n+1} \sum_{i=0}^n {\expe{g}{\pi_{-\lam \Sigma_i}}} g
		= p \ga + (1-p) \gb
		\equiv \tyb + p (\tya-\tyb).
	\endarc
So we have 
	\begar
	\tyb \le \hg_{\aaa} \le \tyb + C \frac{\log n}{n} \le \tya,
	\endar
which is much stronger than \eqref{eq:bbbxprop} (and much simpler to prove).

\iffalse TMP
Besides we recall that for the least square loss, the progressive indirect mixture rule
can be used for $\lam'$ slightly larger than $\lam$ and this leads to \eqref{eq:1}
with a better constant $C$ than the one obtained for the progressive mixture rule 
(see \cite[Example 3.13]{Hau98} and \cite[Section 4]{Aud06a}).

Both assertions go in the same direction: it is likely 
that a progressive indirect mixture rule, e.g. the one in which 
$\hh_i(x)$ maximizes (or almost maximizes) 
	\begar
	y'\mapsto\und{\sup}{y\in\Y} \big\{ 
		- \frac{1}{\lam'} \log {\expe{g}{\pilam}} e^{-\lam' \ell[y,g(x)]} - \ell(y,y') \big\}
	\endar
for a $\lam'$ possibly larger than $\lam$, outperforms significantly the progressive mixture rule on real data.
\fi
\end{remark}

\appendix

\section{Proof of Theorem \ref{th:1}} \label{sec:proof1}

The first assertion is a direct consequence of Lemma 3.3 and Corollary 4.1 of \cite{Aud06a}.
The second assertion is based on an Assouad's type lower bound
(\cite[Inequality (8.19)]{Aud06b}.
Let $\yb=2a-\ya$ and $\tm=\lfloor \log_2 |\G| \rfloor.$
We use the notation introduced in \cite[Section 8.1]{Aud06b}.
We consider a $\big( \tm , \frac{1}{n+1} \wedge \frac{1}{\tm} , 1 \big)$-hypercube
of proba\-bility distributions with $\ha\equiv \argmin_{y\in\Y} \ell_{\ya}(y)$ and $\hb\equiv 
\argmin_{y\in\Y} \ell_{\yb}(y)$. %=2a-\ha$.
We obtain
	\begar
	\E R(\hg) - \und{\min}{g\in\G} R(g) 
		& \ge & \big( \frac{\lfloor \log_2 |\G| \rfloor}{n+1} \wedge 1\big) 
				\da \big(1-\frac{1}{n+1} \wedge \frac{1}{\lfloor \log_2 |\G| \rfloor}\big)^n \\
		& \ge & \big( \frac{\lfloor \log_2 |\G| \rfloor}{n+1} \wedge 1\big) \da e^{-1},
	\endar
where the last inequality comes from $[1-1/(n+1)]^{n} \searrow e^{-1}$. Now the edge discrepancy $\da$ can be computed:
	\begar
	\da & = & \psi_{1,0,\ya,\yb}(1/2) \\
	& = & \und{\inf}{y\in\Y} \frac{ \ell(\ya,y)+\ell(\yb,y) }{2} - \demi
		\und{\inf}{y\in\Y} \ell(\ya,y) - \demi \und{\inf}{y\in\Y} \ell(\yb,y)\\
	& = & \und{\inf}{y\in\Y} \frac{ \ell(\ya,y)+\ell(\ya,2a-y) }{2} - \und{\inf}{y\in\Y} \ell(\ya,y) \\
	& = & \und{\sup}{y\in\Y} [\ell(\ya,a)-\ell(\ya,y)],
	\endar
where the last equality uses that the function $y\mapsto\frac{ \ell(\ya,y)+\ell(\ya,2a-y) }{2}$ is convex.
Finally, from the ``well behaved at center'' assumption, the supremum is positive.

\section{Proof of Theorem \ref{th:2}} \label{sec:proof2}

Let $\tg \in \argmin_{\G} \, R$ and $\eta>0$.
Hoeffding's inequality applied to the random variable $W=\ell[Y,\tg(X)]-\ell[Y,g(X)] \in [-B;B]$
for a fixed $g\in\G$ gives
	\begar
	\E e^{\eta[W-\E W]} \le e^{\eta^2 B^2/2} 
	\endar
for any $\eta>0$. Since the random variable $Z_1,\dots,Z_n$ are independent, we obtain
	\begar
	\E e^{\eta[n R(g) - n R(\tg) + \Sigma_n(\tg) - \Sigma_n(g) ]} 
		& \le & e^{\eta^2 n B^2/2}.
	\endar
Consequently we have
	\begar
	n\big\{ \E R(\hg_{\erm}) - R( \tg ) \big\} & \le & \E \big\{ n R(\hg_{\erm}) - n R(\tg) 
		+ \Sigma_n(\tg) - \Sigma_n(\hg_\erm) \big\}\\
	& \le & \frac{1}{\eta} \log \E e^{ \eta [ n R(\hg_{\erm}) - n R(\tg) 
		+ \Sigma_n(\tg) - \Sigma_n(\hg_\erm) ] }\\
	& \le & \frac{1}{\eta} \log \E \und{\sum}{g\in\G} e^{ \eta [ n R(g) - n R(\tg) 
		+ \Sigma_n(\tg) - \Sigma_n(g) ] }\\
	& \le & \frac{1}{\eta} \log \big( |\G| e^{\eta^2 n B^2/2} \big).
	\endar
The first assertion follows from the (optimal) choice $\eta=\sqrt{\fracb{2\log|\G|}{nB^2}}$.

The second assertion is based on an Assouad's type lower bound. 
Let $\yb=2a-\ya$ and $\tm=\lfloor \log_2 |\G| \rfloor.$
We use the notation introduced in \cite[Section 8.1]{Aud06b}.
We consider a $\big( \tm , \frac{1}{\tm} , \tdb \big)$-hypercube
of probability distributions with $\ha\equiv \tya$ and $\hb\equiv \tyb \eqdef 2a-\tya$ and
$\tdb$ has to be optimized in $[0;1]$.
In the proof of Theorem \ref{th:1}, we take the set $\G$ such that 
$\min_{g\in\G} R(g)=\min_{g} R(g)$, where the second minimum
is $\wrt$ all possible prediction functions. Here the trick is to realize
that $\min_{g\in\G} R(g)$ for our learning setting equals to 
$\min_{g} R(g)$ for the learning task in which the output space is only $\{\tya,\tyb\}$.
Therefore we apply (\cite[Inequality (8.17)]{Aud06b}
with the function $\phi$ appearing in the edge discrepancy $\da$ defined as
	\begarc
	\phi_{y_1,y_2}(p) = \und{\min}{y\in\{\tya,\tyb\}} \big\{ p \ell(y_1,y) + (1-p) \ell(y_2,y) \big\}.
	\endarc
We get
	\begar
	\E R(\hg) & \ge & \und{\min}{g\in\G} R(g) + mw\da\big(1-\sqrt{nw\db}\big)\\
	& = & \und{\min}{g\in\G} R(g) + \da\Big(1-\sqrt{\frac{n}{\tm}\tdb}\Big).
	\endar
From the symmetry and admissibility assumptions of the loss function, we have
	$\ell(\yb,\tyb)=\ell(\ya,\tya) > \ell(\yb,\tya) = \ell(\ya,\tyb)$,
hence
    \begarc
    \delta\eqdef \ell(\ya,\tyb)-\ell(\ya,\tya)>0.
    \endarc
We obtain
	\begar
	\da & = & \psi_{\frac{1+\sqrt{\tdb}}{2},\frac{1-\sqrt{\tdb}}{2},\ya,\yb}(1/2)\\
	& = & \phi_{\ya,\yb}(1/2) - \demi \phi_{\ya,\yb}\Big(\frac{1+\sqrt{\tdb}}{2}\Big)
		- \demi \phi_{\ya,\yb}\Big(\frac{1-\sqrt{\tdb}}{2}\Big)\\
	& = & \phi_{\ya,\yb}(1/2) - \phi_{\ya,\yb}\Big(\frac{1+\sqrt{\tdb}}{2}\Big)\\
	%& = & \ell(\ya,\tya) - \Big( \frac{1+\sqrt{\tdb}}{2} \ell(\ya,\tya) + \frac{1-\sqrt{\tdb}}{2} \ell(\yb,\tya) \Big)\\
	& = & \demi \ell(\ya,\tya)+\demi \ell(\yb,\tya) 
		- \Big( \frac{1+\sqrt{\tdb}}{2} \ell(\ya,\tya) + \frac{1-\sqrt{\tdb}}{2} \ell(\yb,\tya) \Big)\\
	& = & \frac{\sqrt{\tdb}}{2}\delta.
	\endar 
The optimization of the lower bound leads us to choose $\tdb=\frac{\tm}{4n}\wedge 1$ and we get the desired result.

\bibliographystyle{plain}
\bibliography{ref}

\end{document}